
\documentstyle{amsppt}
\pagewidth{6.5in}
\pageheight{9.in}
\topmatter
\title
Asymptotic Stability of Traveling Wave Solutions for Perturbations with Algebraic Decay
\endtitle
\rightheadtext{Perturbations with Algebraic Decay}
\author
Hans Engler 
\endauthor
\affil
Department of Mathematics \\
Georgetown University \\
Washington, D.C. 20057
\endaffil
\email
engler@georgetown.edu
\endemail

\thanks
{Supported by the National Science Foundation (IRD program)} 
\endthanks

\keywords
Generalized Burgers equation, generalized Korteweg-deVries-Burgers equation, traveling wave,  stability, weighted norm
\endkeywords

\subjclass 
35Q53
\endsubjclass

\abstract
For a class of scalar partial differential equations that 
incorporate convection, diffusion, and possibly dispersion 
in one space and one time dimension, the stability of traveling wave 
solutions is investigated.  If the initial perturbation of the traveling 
wave profile decays at an algebraic rate, then the solution is shown to 
converge to a shifted wave profile at a corresponding temporal algebraic 
rate, and optimal intermediate results that combine temporal and spatial 
decay are obtained.  The proofs are based on a general interpolation principle 
which says that algebraic decay results of this form always follow if exponential 
temporal decay holds for perturbation with exponential spatial decay and 
the wave profile is stable for general perturbations.
\endabstract

\refstyle{C}

\endtopmatter

\document

\head
1. Introduction
\endhead

The topic of this note is the class of regularized scalar conservation laws in one spatial dimension 
$$
u_t + \beta u_{xxx} - \alpha u_{xx} + g(u)_x = 0
$$
where subscripts denote partial derivatives. The regularization is due to the presence of viscous terms ($\alpha > 0$) and dispersive terms ($\beta \neq 0$).  The case $g(u) = u^2/2$ is typical and has received much attention.  If $\alpha > 0 = \beta$, this is known as {\smc Burgers} equation.  If  $\alpha = 0 < \beta$, this is essentially the {\smc Korteweg - deVries} equation.  The case $\alpha, \, \beta  > 0$ thus is referred to with a canonical all-Dutch name; it also been studied extensively, as has been the case of general $g$.  I want to study the stability of traveling wave solutions of the form $u(x,t) = \phi(x-ct)$ with respect to perturbations of the initial data, in the cases $\alpha > 0 = \beta$ and $\alpha, \beta > 0$.  Here $c$ is the speed at which the wave profile $\phi$ travels to the right (if $c>0$).  Only monotone wave profiles $\phi$ will be considered. The question is whether the solution approaches a traveling wave in some sense. A natural setting for this is a spatial coordinate system that moves along with the expected wave profile at speed $c$. 

\medskip
Since all shifted wave profiles also give rise to traveling wave solutions, one can only expect that the solution will converge to some shifted profile  $\phi(x-ct-h)$.  Since the quantity 
$\int_\bold R \left( u(x,t) - v(x,t) \right) \, dx$ is independent of $t$ for any pair of solutions $u$ and $v$ for which it is finite, the shift $h$ must be such that $\int_\bold R \left( u(x,0) - \phi(x-h) \right) \, dx = 0$. It is easy to see that the quantity on the left hand side is an affine function of $h$, and thus $h$ can be determined explicitly and can be considered a known quantity.  

\medskip
Stability in this sense was first studied in \cite{6}, for the viscous case $\beta = 0$.  These authors noted that one cannot expect a rate of convergence that holds for all classes of perturbations.  However, in 1976, Sattinger showed in \cite{15} that an exponential rate of decay holds in a moving coordinate system if the perturbation of the initial value decays exponentially in space.  Since this class of perturbations is somewhat restrictive, one may ask the question what the consequences of algebraic decay of the initial perturbation are. It was shown in \cite{8}, \cite{10} and \cite{11} that in this case algebraic decay for the initial perturbation can be "traded in" for some temporal algebraic decay.  Heuristically, the equation for the perturbation behaves like $w_t \mp w_x = 0$ near $x = \pm \infty$ in this situation, that is, the solution near $\infty$ behaves like $w(x,t) \sim w(x+t,0)$.  Therefore, if $w(x,0) \sim e^{-x}$, then $w(x,t) \sim e^{-t} e^{-x}$ near $x = \infty$ as $t \to \infty$ (exponential decay in exponentially weighted norms), and if $w(x,0) \sim x^{-k}$, then $w(x,t) \sim x^{m-k}t^{-m}$ near $x = \infty$ as $t \to \infty$ (algebraic decay with a weaker algebraic weight). For finite $x$, diffusion dominates and leads to exponential decay. Similar results were shown for case $\beta > 0$ in \cite{13}, assuming the wave profile $\phi$ is monotone.  The proofs rely on a detailed study of the spectrum of the linearized problem or on {\it ad hoc} energy estimates. A more refined analysis relying on Green's function estimates is carried out in \cite{5}. 

\medskip
The goal of this note is to show that this "trade-off" follows whenever exponential decay holds in a setting with spatial exponential weights together with (simple) stability in a setting without weights.  Showing these two conditions is often easier than analyzing the full problem in a setting with polynomials weights.  On the other hand, for viscous conservation laws, this interpolation argument can only be applied in the "totally compressive case" where all characteristics run into the shock for the inviscid problem.  This is always true for the scalar case, but systems of regularized conservation laws from physical situations usually do not fall in this category; see \cite{3} and \cite{16} for the much more complicated theory for this case. On the other hand, the argument given in this note is not restricted to the case of one space dimension. 

\medskip
The paper is organized as follows.  In the next section, I show a general interpolation result for linearized problems.  In section 3, this is applied to study the generalized {\smc Burgers} equation. The main result is a sharp decay result in spaces with algebraic spatial decay.  In section 4, a similar result is shown for the generalized {\smc Korteweg-deVries - Burgers} equation.  Appendices A and B contain results for the corresponding linearized problems.  Appendix C contains a simple integral inequality that is used for the passage from linear to nonlinear stability.

\medskip
Here is some notation that is used throughout the paper.  Let $\Omega \subset \bold R^n$ be an unbounded measurable set.  
For $1 \le p < \infty, \, k > 0, \rho > 0$ let us define the function spaces
$$
\align
L^{p,k} &= \{u \in L^p(\Omega) \, \big| \, \int_\Omega |u(x)|^p 
(1+|x|)^{kp} \, dx = \|u\|_{p,k}^p < \infty \} \\
L^{p,\rho} &= \{u \in L^p(\Omega) \, \big| \, \int_\Omega |u(x)|^p 
e^{\rho p |x|} \, dx = \|u\|_{p,\rho}^p < \infty \} 
\endalign
$$
with their natural norms.  The usual modifications are made to define
$L^{\infty,k}$ and $L^{\infty,\rho}$. 
The norms on the unweighted $L^p$-spaces are denoted by $\|\cdot\|_p$. 
The set $\Omega$ is specified to be the real line in sections 3 and 4 and in appendices A and B, thus it does not appear further in the notation.  Constants are denoted by the same letter $C$ whose value may change from line to line, depending only on values that can be expressed in terms of quantities mentioned in the assumptions of a result. If constants have indices, their values remain constant throughout a proof.

\bigskip

\head
2. A Linear Interpolation Result
\endhead

Suppose we are given an operator $S$, not necessarily linear, which maps 
the space $L^p(\Omega)$ into itself and also $L^{p,\rho}(\Omega)$ into itself, with 
estimates
$$
\align
\|S(u) - S(v) \|_p &\le C_0 \|u - v \|_p \qquad \text {for all} \qquad
u, \, v \in L^p(\Omega) \tag 2.1a \\
\|S(u)\|_{p,\rho} &\le C_0 e^{-t} \|u  \|_{p,\rho} \qquad \text 
{for all} \qquad u \in L^{p,\rho}(\Omega) \, \tag 2.1b
\endalign
$$
for some constants  $\rho, \, C_0, \, t > 0$.  
The main result of this section says that $S$ also maps $L^{p,k}$ into
$L^{p,l}$ for $0 < l < k$ and gives an estimate for this mapping.  

\medskip
I shall prove this in detail for the case $p < \infty$.
A different proof will be given for the case where $S$ is linear and
$p = \infty$. 

\medskip
\proclaim{Theorem 2.1}
Let $1 \le p < \infty$. 
Under the above conditions, $S$ maps $L^{p,k}$ into $L^{p,l}$ for all
$0< l < k$, and there exists $C_1 > 0$, 
depending only on $p$, such that for all $v \in L^{p,k}$
$$
 \|Sv\|_{p,l} \le C_0C_1 \rho^{k-l}(1+t)^{l-k} 
\|v\|_{p,k} \, . \tag 2.2
$$
\endproclaim

\medskip
\demo{Proof}
Let us assume $\rho = 1$ and continue to 
write $\| \cdot\|_{p,\rho}$ for the corresponding norm with exponential
weight. Fix $p \in [1, \infty)$ and define for $r \ge 0$
$$
m_p(r) = \frac{r}{\big(1 + r^{\frac{p}{p-1}}\big)^{\frac{p-1}{p}}}
$$
if $p > 1$ and $m_1(r) = \min(1,r)$. 
Clearly $m_1(r) \le 2 m_p(r) \le 2m_1(r)$ for all $r \ge 0$ and all $p$.
For $s \in \bold R$ and $u \in L^p(\Omega)$ define the functional 
$$
K(s,u) = \big(\int_\Omega (|u(x)| m_p(e^{s+|x|}))^p \, dx \big)^{1/p} \, .
$$
This is clearly an equivalent norm on $L^p$ and in fact a modified 
K-functional (\cite{1}), namely
$$
K(s,u) = \inf_{v \in L^{p,\rho}} \big(\|u - v\|_p^p + e^{sp}\|v\|_{p,\rho}^p
\big)^{1/p} \, . \tag 2.3
$$
Indeed,
$$
\align
\inf_{v \in L^{p,\rho}} \big(\|u - v\|_p^p + e^{sp}\|v\|_{p,\rho}^p
\big) &=  \inf_{v \in L^{p,\rho}} \int_\Omega
\big(|u(x) - v(x)|^p + e^{sp + |x|p}|v(x)|^p \big) \, dx \\
 &\ge  \int_\Omega \inf_\zeta
\big(|u(x) - \zeta)|^p + e^{sp + |x|p}|\zeta|^p \big) \, dx \\
&=  \int_\Omega \big(|u(x)| m_p(e^s+|x|)\big)^p \, dx \\
&=  \big(\|u - v_0\|_p^p + e^{sp}\|v_0\|_{p,\rho}^p\big) \, .
\endalign
$$
Here 
$$
v_0(x) = \frac{u(x)}{1 + e^{(s + |x|)\frac{p}{p-1}}}
$$
for $ p>1$. For $p = 1$ one sets $v_0(x) = u(x) $ for $s + |x| \le 0$
and $v_0(x) = 0$ otherwise, and the last equation is again valid.  
Then $v_0 \in L^{p,\rho}$, and (2.3) follows.
The definition shows immediately that $K(\cdot,u)$ is 
differentiable, non-decreasing, and bounded above by $\|u\|_p$ for all $p$.  
Also, 
$$
2^{-p}\int_\Omega |u(x)|^p m_1(e^{s+|x|})^p \, dx \le K(s,u)^p \le
\int_\Omega |u(x)|^p m_1(e^{s+|x|})^p \, dx \, . 
$$
Next fix also $k>0$ and set
$$
h_k(s) = \cases 	&e^{-s} \quad \quad (s \ge 0) \\
			&(1-s)^{kp-1} \quad \quad (s < 0)
	 \endcases \, . \tag 2.4
$$
and
$$
\|u\|_*^p = \int_{-\infty}^{\infty} K(s,u)^p h_k(s) \, ds 
$$
whenever this quantity is finite. The next claim is that
$$
L^{p,k} = \big\{u \, \big| \, \|u\|_*^p < \infty \, \big\} 
$$
and that $\| \cdot \|_*$ is an equivalent norm on this space.
Indeed,
$$
\align
\int_{-\infty}^{\infty} K(s,u)^p h_k(s) \, ds &\le
\int_{-\infty}^{\infty} \int_\Omega |u(x)|^p \min(1, e^{sp + |x|p})
h_k(s) \, dx \, ds\\
& = \int_\Omega |u(x)|^p \left(\int_{-\infty}^{-|x|}e^{sp+|x|p}
(1-s)^{kp-1} \, ds + \int_{-|x|}^0(1-s)^{kp-1} \, ds + 
\int_0^{\infty} e^{-sp} \,ds \right) \, dx \\
&\le  \int_\Omega |u(x)|^p C\left(1 + |x|^{kp} \right) \, dx \le
C \|u\|_{p,k}^p \, .
\endalign
$$
Reversely,  
$$
\align
\int_{-\infty}^{\infty} K(s,u)^p h_k(s) \, ds &\ge
2^{-p} \int_{-\infty}^{\infty} \int_\Omega |u(x)|^p \min(1, e^{sp + |x|p})
h_k(s) \, dx \, ds\\
& \ge 2^{-p} \int_\Omega |u(x)|^p \int_{-|x|}^0(1-s)^{kp-1} \, ds  
\, dx \\
&\ge  C \|u\|_{p,k}^p \, .
\endalign
$$
Let now  $S$ be an operator satisfying (2.1a,b).  Then
$$
\align
K(s,Su)^p &= \inf_{v \in L^{p,\rho}} \left(\|Su - v\|_p^p +
e^s\|v\|_{p,\rho}^p \right) \\
&\le \inf_{v \in L^{p,\rho}} \left(\|Su - Sv\|_p^p +e^s\|Sv\|_{p,\rho}^p \right)\\
&\le C_0\inf_{v \in L^{p,\rho}} \left(\|u - v\|_p^p +e^{s-t}\|v\|_{p,\rho}^p \right) \\
&= C_0 K(s-t,u)^p \, .
\endalign
$$
Let $0 < l < k$, let $u \in L^{p,\rho}$, and set $H_l(r) = \int_r^{\infty} h_l(\tau) \, d \tau$ and 
$k(s) = \frac{d}{ds} K(s,u)^p \ge 0$.  Then 
$$
\align
\|Su\|_{p,l}^p &\le C\int_{-\infty}^\infty K(s,Su)^p h_l(s) \, ds \\ 
  & \le CC_0\int_{-\infty}^\infty K(s-t,u)^p h_l(s) \, ds \\
  & = CC_0\int_{-\infty}^\infty k(s) H_l(s+t) \, ds  \, .
\endalign
$$
An elementary calculation shows that
$H_l(s+t) \le C H_k(s)(1+t)^{l-k}$ for all $s$ and $t$. One can therefore estimate further
$$
\|Su\|_{p,l}^p  \le CC_0(1+t)^{l-k}\int_{-\infty}^\infty k(s) H_k(s) \, ds  = 
C_1 C_0 (1+t)^{l-k} \|u\|_{p,k} \, .
$$ 
This proves the theorem in the case $\rho = 1$.  The general case follows by a 
scaling argument. 
\enddemo

\medskip
The proof can be modified to extend to the case $p = \infty$. However, I 
prefer to give an alternative proof in this case.  It extends to 
subspaces of $L^\infty$ that are closed under multiplication with 
smooth functions that are bounded together with their derivatives.  
Such spaces include $X = BC^m(\Omega)$ and 
$X = L^\infty(\Omega) \cap UC(\Omega)$, where $BC^m$ is the space of
$m$-times differentiable functions with bounded derivatives and $UC$
is the set of uniformly continuous functions on $\Omega$. 
The result is formulated for the cases $X = L^\infty(\Omega)$, $X = BC^0(\Omega)$, 
and $X = L^\infty(\Omega) \cap UC(\Omega)$, equipped with the supremum
norm. Accordingly let $X_k = X \cap L^{\infty,k}$ and $X_\rho = X \cap
L^{\infty,\rho}$, equipped with their natural norms 
$\|\cdot\|_{\infty,k}$ and $\|\cdot\|_{\infty,\rho}$.  

\medskip
\proclaim{Theorem 2.2} 
Let $S:X \, \to \, X$, $S:X_\rho \, \to \, X_\rho$ be a linear operator 
for which (2.1a,b) holds. Then $S$ maps $X_k$ into $X_l$ for all
$0< l < k$, and there exists a constant $C_1 > 0$ depending on $k$ and $\rho$
such that for all $v \in X_k$
$$
 \|Sv\|_{\infty,l} \le \frac{C_0C_1}{(1+t)^{k-l}} 
\|v\|_{\infty,k}  \, . \tag 2.5
$$
\endproclaim

\medskip
\demo{Proof}
Let $v \in L^{\infty,k}$ be given. Thus $(1+|x|)^k|v(x)| \le A$ almost
everywhere for some smallest constant  $A$.  Let $w = Sv$.  
The goal is to show that 
$(1+t)^{k-l}(1+|x|)^k|w(x)| \le C_0C_1A$ almost everywhere for some constant $C_1$.  
Let $R \ge 0$, to be chosen later.  
Choose $\varphi \in C^\infty(\bold R)$ with $\varphi(r) = 0$ 
for $r\le 0$, $\varphi(r) = 1$ for $r \ge 1$, and $\varphi' \ge 0$. Set 
$v_1(x) = v(x) \varphi(|x| - R)$ and $v_2 = v-v_1$. Then $\|v_1\|_\infty \le 
A(1+R)^{-k}$ and $\|v_2\|_{\infty,\rho} \le C_0C(1+R)^{-k}e^{\rho R}$ with
$C = C(k,\rho) \ge 0$. Thus
$$
\|S v_1 \|_\infty \le C_0 (1+R)^{-k}  A \quad \text{and} \quad
\|S v_2 \|_{\infty,\rho} \le C_0C(1+R)^{-k}e^{\rho R-t} A \, .
$$
Let $x \in \Omega$.  Then
$$
\align
(1+t)^{k-l}(1+|x|)^l |w(x)| & \le (1+t)^{k-l}(1+|x|)^l\left(|Sv_1(x)| + 
|Sv_2(x)| \right)\\
	&\le C_0(1+t)^{k-l}(1+|x|)^l(1+R)^{-k}A(1+e^{\rho R - \rho|x| - t}) \, .
\endalign
$$
Set $\sigma = \min \left( \dfrac{l}{k}, \, \dfrac{k-l}{\rho k} \right)$
and choose $R = \sigma t^\frac{k-l}{k} |x|^{\frac{l}{k}} \le \frac{t}{\rho}
 + |x|$. One can estimate further
$$
(1+t)^{k-l}(1+|x|)^l |w(x)| \le  2 C_0(1+t)^{k-l}(1+|x|)^l (1+R)^{-k} A 
	\le C_0C_1 A \, .
$$
This proves the theorem.
\enddemo

\medskip
At first glance, it is surprising that the algebraic decay estimates are 
independent of the exponential spatial weight that appears in the assumptions. 
The scaling argument used at the end of the proof explains this phenomenon and shows where the constant $\rho$ reappears in the result.
The special case $\Omega = (-\infty,0]$ with the right shift operator $S_tv(x) = v(x-t)$ for  
$x \le 0$ shows that all estimates in the two theorems are sharp, up to the values of the constants.  Indeed, the use of K-functionals in the proof of theorem 2.1 makes the argument resemble the direct proof for this special case.

\bigskip

\head
3. Scalar Viscous Conservation Laws
\endhead

In this section scalar viscous conservation laws of the following form are considered: 
$$
u_t - u_{xx} + f(u)_x = 0 \quad (x \in \bold R, \, t > 0)\, .\tag 3.1
$$
Here $f:\bold R \, \to \, \bold R$ is a $C^2$-function with uniformly
bounded derivatives.  A traveling wave solution $u(x,t) = \phi(x-ct)$ with
speed $c$ and limiting behavior $\phi(r) \, \to \, \phi_\pm$ as $r \,  \to \, \pm \infty$ is easily seen to exist if and only if $c$ equals the slope of
the line segment connecting the points $(\phi_+, f(\phi_+))$ and
$(\phi_-,  f(\phi_-))$ and the graph of $f$ lies entirely above or below 
this line segment.  Accordingly, the wave profile $\phi$ is decreasing 
or increasing. After rescaling the dependent variable $u$, adding a linear
function to $f$, and changing to
a moving coordinate system, one can assume that $c = 0, \, \phi_- = 1, \, 
\phi_+ = 0$. The prototypical example is $f(r) = r^2 - r$, with the wave
profile $\phi(r) = \left( 1+e^r \right)^{-1}$. 

\medskip
The focus is the convergence of solutions
of (3.1) to some translate $\phi(\cdot - h)$ as $t \, \to \, \infty$, where $h$ is known.
Let us therefore assume that $h = 0$.
It was shown in \cite{6} that this convergence in the uniform sense follows
if $f$ is uniformly convex and $\int_0^\infty |u(x,0)| \, dx +
\int_{-\infty}^0 |u(x,0)-1| \, dx$ is finite.  If $f$ is merely $C^2$-smooth but not necessarily convex, the same conclusion holds in the $L^1$-sense even for the case of quasilinear diffusion (\cite{13}). In the seminal paper
\cite{15}, it was shown that convergence at an exponential rate holds in spaces
with exponentially weighted norms, namely
$$
\|u(\cdot,t) - \phi\|_{\infty, \epsilon} = O(e^{-\delta t}) \tag 3.2
$$
for some $\epsilon, \, \delta > 0$, if the quantity on the
left hand side is sufficiently small for $t = 0$. This holds for arbitrary
non-convex $f$, assuming only that 
$$
f'(0) \ne 0 \ne f'(1) \, . \tag 3.3
$$
In several recent papers, the stability of wave profiles in spaces with
polynomial weights was discussed.  
Assuming only (3.2) it was shown in \cite{8} that if
$$
\|u(\cdot,0) - \phi\|_{\infty,k+m} = \delta \tag 3.4a
$$
is sufficiently small, then 
$$
\|u(\cdot,t) - \phi(\cdot)\|_{\infty,k} \le C (1+t)^{-m/2} \delta \tag 3.4b
$$
for integers $k,m$ satisfying $k \ge 1, \, 2 \le m \le k+1$, or $k = 1, \,
m \ge 2$. A comparable result for this situation from \cite{10} and \cite{11}
assumes that
$$
\|\Psi\|_{2,\alpha} < \infty \tag 3.5a
$$
where $\Psi(x) = \int_{-\infty}^x \left( u(z,0) - \phi(z) \right) \, dz$,  
and that a suitable unweighted $L^2$ - norm of $u(\cdot,0) - \phi$
is small. The conclusion then is essentially that
$$
\|u(\cdot,t) - \phi\|_\infty \le C (1+t)^{-\alpha}  \, . \tag 3.5
$$
Here $\alpha > 0$ is arbitrary. Both results show a trade-off between
the spatial decay of the initial data and the temporal decay of the
solution. A much more detailed and general result in \cite{5} implies 
that spatial decay of the antiderivative $\Psi$ dominates the 
temporal decay of $u(\cdot,t) - \phi$ for finite $x$, while the spatial
decay of $u(\cdot,0) - \phi$ dominates the temporal and spatial decay as 
both $x$ and $t$ go to $\infty$, with the canonical trade-off.  

\medskip
The main result for this situation assumes also (3.3).

\medskip
\proclaim{Theorem 3.1} Let $k>1, \, 0 < m < k$ be real numbers, and set 
$\Psi(x) = \int_{-\infty}^x \left( u(z,0) - \phi(z) \right) \, dz$.
There exists a constant $C_0$ such that if
$$
\|\Psi\|_{\infty,k}  = \epsilon (u_0) \tag 3.6a
$$ 
is sufficiently small, then the solution $u$ exists for all $t>0$, and for all $0 < m < k$
$$
\|u(\cdot,t) - \phi\|_{\infty,m} \le C_0  (1+t)^{m-k}
\epsilon(u_0)  \, . \tag 3.6b
$$ 
\endproclaim

\medskip
\demo{Proof}
The proof follows a pattern which will be repeated in the next section.  
A formal linearization of the problem is introduced (step 0),
a function space setting is defined, and a solution $u$ is produced, 
using results about the linearized equation from Appendix A (step 1).
After specifying the short time behavior of the solution (step 2),  
suitable {\it a priori} estimates are shown with the aid of Lemma C.1, 
which finishes the proof. 
 
\medskip
{\bf Step 0.} Let $u$ be a solution
of (3.1), and set $v(x,t) = \int_{-\infty}^x \left( u(z,t) - \phi(z) \right)
\, dz$. Thus $v(x,0) = \Psi(x)$, and $v$ solves the equation
$$
v_t(x,t) - v_{xx}(x,t) = f(\phi(x)) - f(\phi(x) + v_x(x,t)) \tag 3.7
$$
or
$$
v_t(\cdot,t) + Lv(\cdot,t) = F(v)(\cdot,t)
$$
where formally
$$
\align
Ly(x) &= -y_{xx}(x) + f'(\phi(x))y_x(x)  \tag 3.8 \\
F(y)(x) &= f(\phi(x)) + f'(\phi(x))y_x(x)- f(\phi(x) + y_x(x)
= f''(\zeta(x))y_x^2(x) \, . \tag 3.9
\endalign
$$
Here $\zeta(x)$ is a number between $\phi(x)$ and $\phi(x) + y_x(x)$.

\medskip
{\bf Step 1.} Recall the notation and the results from Appendix A:
$X = L^\infty(\bold R) \cap C^0(\bold R)$, $X_0 = \{y \in X \, \big| \,
y \, \text{is uniformly continuous} \}$ are Banach spaces, both equipped with the 
supremum norm. The operator $L$ acting in $X$ is defined in (3.8), and the restriction 
of $L$ to $X_0$ is also denoted by $L$. Then $-L$ generates an analytic semigroup 
$(S(t))_{t \ge 0}$ in  $X_0$ and in various weighted spaces, especially in 
the spaces $X_\rho$ and $X_k$. 
From the properties of $f$ and the definition of $F$ in (3.8), 
one immediately deduces the estimates
$$
\align
\|F(y_1)- F(y_2) \|_{\infty,k} &\le C\|y_{1,x} - y_{2,x} \|_{\infty,k} \tag 3.10a \\
\|F(y)\|_{\infty,k} &\le C\|y_x \|_{\infty,k/2}^2  \tag 3.10b
\endalign
$$
for any $k \ge 0$ and suitable $y, \, y_1, \, y_2$. Suppose that $\Psi \in 
X_k \subset X_0$ is given and satisfies the assumptions of Theorem 3.1. 
Let us seek a continuous $X_0$-valued
solution $v(\cdot,t)$ of the integral equation
$$
v(\cdot,t) = S(t)\Psi + \int_0^t S(t-s) F(v)(\cdot,s) \, ds \, .\tag 3.11
$$
A contraction argument, using (3.10a), produces such a solution on some 
finite time interval, and the solution satisfies (3.7) in the classical 
sense if it is sufficiently smooth. Due to the global Lipschitz property 
in (3.10a), this solution exists for all times. The existence and 
uniqueness arguments hold in all spaces $X_k$, $0 < k \le m$, 
and therefore $v(\cdot,t)$ belongs to all these spaces for all $t$.  
Set $w(x,t) = u(x,t)-\phi(x) = v_x(x,t)$, and let us write 
$N(t) = \max (1, t^{-1/2} )$. \medskip

{\bf Step 2.}
Let us now characterize the behavior of $w$ for 
$0<t \le 1$ in more detail. The goal is to show the estimate
$$
\|w(\cdot,t)\|_{\infty,m} \le CN(t) \epsilon(u_0) \quad (0 < t \le 1) \tag 3.12
$$
for all $0 < m \le k$.
For this purpose consider the integral
equation for $w$ which is obtained by differentiating (3.1), namely
$$
w(\cdot,t) = [S(t)\Psi]_x + \int_0^t [S(t-s) F(v)(\cdot,s)]_x \, ds \, .\tag 3.13
$$
For the ``free term'' $[S(t)\Psi]_x$ in (3.13), (3.12) is just estimate (A.6). 
Using (3.12), (A.6), and (3.10b) one then derives the inequality
$$
\frac{\|w(\cdot,t)\|_{\infty,m}}{N(t)} \le C\epsilon(u_0)  + 
\int_0^t C \frac{N(t-s)N(s)}{N(t)} \frac{\|w(\cdot,s)\|_{\infty,m}}{N(s)} \, ds \,. \tag 3.14
$$
A standard argument for linear integral inequalities implies that $\frac{\|w(\cdot,t)\|_{\infty,m}}{N(t)} \le C$ on $(0,1]$, i.e. (3.12). 

\medskip
{\bf Step 3.} Finally, let us prove estimate (3.6b). The first thing to notice is again
that the estimate
(3.6b) in Theorem 3.1 holds for the ``free'' term $[S(t) \Psi]_x$,
since $\|[S(t)\Psi]_x\|_{\infty,m} \le C \|S(t-1/2)\Psi\|_{\infty,m} \le C(1+t)^{m-k}\epsilon(u_0) $.
Consider first the special case $m = k/2$ and thus $m-k=-m$. Define the quantity
$$\gamma(t) = \sup_{1 \le s \le t} s^m\|w(\cdot,s)\|_{\infty,m}  \tag 3.15
$$
for $t \ge 1$.
Then $\gamma$ is continuous, $\gamma(1) = \|w(\cdot,1)\|_m \le C_0\epsilon(u_0)$ 
for some constant  $C_0$, and 
$\|F(v)(\cdot,s)\|_{\infty,m} \le C s^{-k} \gamma^2(s) $ for all $s$.  Then (3.13)
implies for $t \ge 1$
$$
\align
\|w(\cdot,t)\|_{\infty,m} &\le C t^{-m}\epsilon(u_0) + \int_0^1 C (1+t-s)^{-m} N(t-s)
\|w(\cdot,s)\|_{\infty,m} \, ds  \\ & \quad + \int_1^t C (1+t-s)^{-m}  
N(t-s)s^{-k} \gamma^2(s) \,ds \\  
&\le C_1 t^{-m}\epsilon(u_0) + C_2 t^{-m} \gamma^2(t) 
\endalign
$$
by Lemma C.1 and estimate (A.6). Therefore $\gamma(t) \le C_1\epsilon(u_0) + C_2\gamma^2(t)$
for all $t \ge 1$. If $4C_1C_2\epsilon(u_0)< 1$ and $2C_0C_2\epsilon(u_0) < 1$ 
(thus $\gamma(1) < \left(2C_2\right)^{-1}$),  then $\gamma(t)\le C_0 
\epsilon(u_0)$ for all $t$ by an elementary algebra argument.  This is the desired 
estimate for $m =k/2$, and it holds if $\epsilon(u_0)$ is sufficiently small.
If $0<m<k$ is arbitrary, then one can use the estimates for $w$ on $(0,1]$ and for 
$[S(t)\Psi]_x$ to obtain for $t \ge 1$ 
$$
\align
\|w(\cdot,t)\|_{\infty,m} &\le Ct^{m-k} \epsilon(u_0) + \int_0^1 C (1+t-s)^{m-k} N(t-s)
\|w(\cdot,s)\|_{\infty,m} \, ds \\
& \quad + \int_1^t C(1+t-s)^{m-k} N(t-s) \|F(v)(\cdot,s)\|_{\infty,k} \, ds \\
&\le Ct^{m-k} \epsilon(u_0) + \int_1^t C (1+t-s)^{m-k} N(t-s) \|w\cdot,s)\|_{\infty,k/2}^2 \, ds \\
&\le  Ct^{m-k} \epsilon(u_0) + C \epsilon(u_0)^2 \int_1^t (1+t-s)^{m-k} N(t-s) (1+s)^{-k} \, ds \\
& \le Ct^{m-k} \epsilon(u_0)
\endalign
$$
by the estimate for  $\|w\cdot,s)\|_{\infty,k/2}$ that was just established and by Lemma C.1. The theorem is now completely proved. 
\enddemo

\bigskip
\head
4. Generalized Korteweg-DeVries - Burgers Equations
\endhead

Let us now look at the partial differential equation
$$
u_t - \alpha u_{xx} +u_{xxx} + g(u)_x = 0 \quad (x \in \bold R, \, t > 0)\, .\tag 4.1
$$
The parameter $\alpha$ is positive, and $g$ is $C^2$-smooth.  The case 
$g(u) = (p+1)^{-1}u^{p+1}$ with integer $p > 0$ is a model for long wave 
propagation in media with dissipation and dispersion.  The special case  $p=1$ is
known as {\smc Korteweg-DeVries - Burgers } equation.  It reduces to the 
{\smc Korteweg-DeVries} equation if $\alpha = 0$.  
Under certain conditions, the equation admits monotone traveling wave solutions 
$u(x,t) = \phi(x-ct)$ with
speed $c$ that connect the end states $\phi_\pm = \lim_{r \to \pm \infty} \phi(r)$.
Such a wave profile must satisfy the third order ordinary differential equation
$$
-c \phi' + g(\phi)' + \phi''' - \alpha \phi'' = 0 \, . \tag 4.2
$$
An example is $g(r) = 2r(r-1)(b-r)$ with $b \ge 2$, which has the wave profile
$\phi(r)= \left(1+e^r\right)^{-1}$ for the parameter $\alpha = 2b-1$ and the speed $c = 0$. 
General profiles (not necessarily monotone) have been constructed in \cite{2} and \cite{7}.  
It is known that monotone profiles exist for $g(u) = (p+1)^{-1}u^{p+1}$ and 
$\alpha \ge 2\sqrt{pc}$.  A slightly more situation is the setting for the next result.

\medskip
\proclaim{Proposition 4.1}
Let $g \in C^2$ be strictly convex.  
A monotone wave profile $\phi$ for (4.1) exists if and only if
$$
\align
c &= \frac{g(\phi_+) - g(\phi_-)}{\phi_+ - \phi_-}  \tag 4.3a \\ 
\alpha &\ge 2\sqrt{g'(\phi_-)-c} \tag 4.3b \\
\phi_+ &< \phi_-  \, .\tag 4.3c
\endalign
$$
The profile $\phi$ must therefore be monotonically decreasing.
\endproclaim

\medskip

\demo{Proof} Suppose $\phi$ is a monotone wave profile with limits  $\phi_{\pm}$ at
$r = \pm \infty$.  Clearly  $-c\phi + g(\phi) + \phi'' - \alpha \phi' = const.$ and thus 
$-c\phi_- + g(\phi_-) = -c\phi_+ + g(\phi_+)$, implying (4.3a).  
Set $\psi(z) = \phi_- - \phi(-z)$ and 
$$
f(r) = g(\phi_-) - c r -g(\phi_- -r)\, . \tag 4.4
$$
Then $f$ is concave, and $-\alpha \psi' - \psi'' = f(\psi)$.  This is the equation 
for a wave profile $\psi$ of the {\smc Fisher - Kolmogorov-Petrovskii-Piskunov (F-KPP)} 
equation $v_t - v_{xx} = f(v)$ that
travels to the right with speed $\alpha$ and has limits 
$\psi_- = \phi_- - \phi_+, \, \psi_+ = 0$. It is known (\cite{4}) that such a monotone 
wave profile for concave $f$ exists if and only if $\alpha \ge 2 \sqrt{f'(0)} = 
2 \sqrt{g'(\phi_-) -c}$.  In this case, $\psi_- > \psi_+$ and therefore $\phi_- > \phi_+$,
since $f$ is positive between $\psi_-$ and $\psi_+$.  Thus (4.3b) and (4.3c) are true. 
Conversely, let (4.3a-c) hold.  Define $f$ as in (4.4).  By well-known results about the
{\smc F-KPP} equation, there exists a unique decreasing wave profile $\psi$ with $\psi(0)
= (\phi_- - \phi_+)/2$ that moves to the right with speed $\alpha$.  Then $\phi(z) = \phi_- -
\psi(-z)$ is a monotone wave profile for (4.1) with $\phi(\pm \infty) = \phi_{\pm}$.
\enddemo

\medskip
It is easy to see that in fact  $\phi' < 0$ on $\bold R$. 
After rescaling the independent variable, adding a linear function to  $g$, 
and changing to a moving coordinate system, one can assume that $c=0, \, \phi_- = 1$, 
and $\phi_+ = 0$. As in section 3, the focus is on the convergence of solutions
of (4.1) to some translate $\phi(\cdot - h)$ as $t \, \to \, \infty$. 
As before, we can assume that $h = 0$ and define $\Psi(x) = \int_{-\infty}^x
\left(u(z,0) - \phi(z) \right)\, dz$. 
In \cite{2}, (4.1) was discussed in the case $g(x) = x^2$, and it was shown that 
the difference $u(\cdot,t) - \phi$ converges to 0 in $L^2$ together with 
its derivatives if this difference 
is small in $L^{2,k}$ for some $k>1$ at $t=0$ and if sufficiently many derivatives of 
its derivatives are small in $L^2$. One of the results in \cite{14} says that this
convergence in fact holds if the initial difference is small just in $L^2$. 
The main result in \cite{12} also covers the case $g(u) = u^2$ and states essentially 
that if $\|\Psi(\cdot)\|_{2,k}$ is sufficiently small, then  $\|u(\cdot,t) - \phi\|_{2,m}
= O(t^{m-k+\epsilon})$.  Here $\epsilon = 0$ if $2m-2k$ is an integer, and it is positive but arbitrarily small otherwise.  Derivatives of $u(\cdot,t) - \phi$ are shown to decay at higher rates.

\medskip

The main result of this section assumes that
$$
g'(0) < 0 < g'(1) \quad \text{and} \quad g''(r) > 0 \quad \text{for all} \quad x \, .\tag 4.5
$$

\medskip
\proclaim{Theorem 4.2} Let $k>1$ be a real number.
There exists a constant $C_0$ such that if
$$
\|\Psi\|_{2,k} = \epsilon (u_0) \tag 4.6a
$$ 
is sufficiently small, then the solution $u$ exists for all $t>0$, and for all $0 < m < k$ and $t \ge 1$
$$
\align
\|u(\cdot,t) - \phi\|_{4,m} &\le C(1+t)^{m-k}\epsilon(u_0) \tag 4.7a \\
\|u(\cdot,t) - \phi\|_{2,m} &\le C(1+t)^{m-k}\epsilon(u_0) \tag 4.7b
\endalign
$$ 
\endproclaim

\medskip
\demo{Proof}
The proof follows the scheme used in section 3.    
Two different function space settings ($L^4$ and $L^2$) are used to handle the quadratic
nonlinearity.  Estimates  (B.7a, b) connect these settings.
  
\medskip
{\bf Step 0.} Let $u$ be a solution
of (4.1), and set $v(x,t) = \int_{-\infty}^x \left( u(z,t) - \phi(z) \right)
\, dz$. Thus $v(x,0) = \Psi(x)$, and $v$ satisfies
$$
v_t(\cdot,t) + Lv(\cdot,t) = G(v)(\cdot,t)
$$
where now formally
$$
\align
Ly(x) &= -\alpha y_{xx}(x) + y_{xxx}(x,t) + g'(\phi(x))y_x(x)  \tag 4.8 \\
G(y)(x) &= g(\phi(x)) + g'(\phi(x))y_x(x)- g(\phi(x) + y_x(x)
= g''(\zeta(x))y_x^2(x) \, . \tag 4.9
\endalign
$$

\medskip
{\bf Step 1.} The operator  $L$ acting in $L^2(\bold R)$ is defined in (4.8), and
$-L$ generates a $C_0$ semigroup $S(t)_{t \ge 0}$ in this space by Appendix B. 
Let us note that 
$\frac{d}{dx} g'(\phi(x)) = g''(\phi(x)) \phi'(x) < 0$, and thus $S$ is a contraction semigroup by (B.10). The semigroup can be restricted to the weighted spaces $L^{2,k}$ and $L^{2,\rho}$ ($\rho < \alpha/3$).  From the properties of $g$ and the definition of $G$ in (4.8), one immediately deduces the estimates
$$
\align
\|G(y_1)- G(y_2) \|_{2,k} &\le C\|y_{1,x} - y_{2,x} \|_{2,k} \tag 4.10a \\
\|G(y)\|_{2,k} &\le C\|y_x \|_{4, k/2}^2  \tag 4.10b
\endalign
$$
for any $k \ge 0$ and suitable $y, \, y_1, \, y_2$, with some universal constant  $C$.  
Suppose now that $\Psi \in 
L^{2,k}$ is given. A continuous $L^2$-valued solution $v(\cdot,t)$ of the integral equation
$$
v(\cdot,t) = S(t)\Psi + \int_0^t S(t-s) G(v)(\cdot,s) \, ds \, \tag 4.11
$$
is again found by a contraction argument, using (4.10a), and due to the global Lipschitz property in (4.10a), this solution exists for all times. The existence and 
uniqueness arguments hold in all spaces $L^{2,m}$, $0 < m \le k$, and therefore $v(\cdot,t)$ belongs to all these spaces for all $t$. The solution satisfies (4.7) in the classical 
sense if it is sufficiently smooth. Although the linear part of (4.7) does not enjoy maximal regularity properties, more smoothness for the solution follows easily from smoothness of the data, if the equation is differentiated and the results in Appendix B are used.    
Set $w(x,t) = u(x,t)-\phi(x) = v_x(x,t)$ and write $N_0(t) = \max (1, t^{-1/2} )$ and
$N_1(t) = \max (1, t^{-5/8} )$. \medskip

{\bf Step 2.}
As before, let us next characterize the behavior of $w$ for 
$0<t \le 1$ in more detail. The goal is to show the estimates
$$
\align
\|w(\cdot,t)\|_{2,m} \le CN_0(t) \epsilon(u_0) \tag 4.12a \\
\|w(\cdot,t)\|_{4,m} \le CN_1(t) \epsilon(u_0) \tag 4.12b
\endalign
$$
for $0 < t \le 1$. For this purpose consider the integral
equation for $w$ which is obtained by differentiating (4.1), i.e. (3.13) with $F$ replaced by $G$. 
For the ``free term'' $[S(t)\Psi]_x$ in (3.13), (4.12a, b) follows directly from (B.7a, b). 
To show (4.12a), one uses (B.7a) and (3.10b) to derive inequality (3.14) with $\|w(\cdot,t)\|_{\infty,m}/N(t)$ replaced everywhere by $\|w(\cdot,t)\|_{2,m}/N_0(t)$.
For (4.12b), one uses (B.7b) and arrive at an inequality like (3.14) for $\|w(\cdot,t)\|_{4,m}/N_1(t)$. 
Standard arguments for linear integral inequalities imply that 
$\|w(\cdot,t)\|_{2,m}/N_0(t) + \|w(\cdot,t)\|_{4,m}/N_1(t) \le C$ on $(0,1]$, i.e. (4.12a, b). 

\medskip
{\bf Step 3.} Finally, let us prove estimate (4.6b). Start by observing that (B.14a,b) are just  estimates (4.6a,b) for the ``free'' term $[S(t)\Psi]_x$. Consider again first the special case  $m = k/2$ and thus $m-k=-m$. 
Define the quantity
$$\gamma(t) = \sup_{1 \le s \le t} s^m\|w(\cdot,s)\|_{4,m}  \tag 4.13
$$
for $t \ge 1$.
Then $\gamma$ is continuous, $\gamma(1) = \|w(\cdot,1)\|_{4,m} \le C_0\epsilon(u_0)$ 
for some $C_0$, and 
$\|G(v)(\cdot,s)\|_{2,m} \le C s^{-k} \gamma^2(s) $ for all $s$.  Equation (4.11), estimate (4.12) and Lemma C.1 imply that for $t \ge 1$
$$
\align
\|w(\cdot,t)\|_{4,m} &\le C t^{-m}\epsilon(u_0) + \int_0^1 C (1+t-s)^{-m} N_1(t-s)
\|w(\cdot,s)\|_{2,m} \, ds  \\ & \quad + \int_1^t C (1+t-s)^{-m}  
N_1(t-s)s^{-k} \gamma^2(s) \,ds \\  
&\le C_1 t^{-m}\epsilon(u_0) + C_2 t^{-m} \gamma^2(t) \, .
\endalign
$$
Therefore $\gamma(t) \le C_1\epsilon(u_0) + C_2\gamma^2(t)$
for all $t \ge 1$. As in section 3,  $\gamma(t)\le C_0 \epsilon(u_0)$ follows if 
$\epsilon(u_0)$ is sufficiently small, which is the desired estimate.
If $0<m<k$ is arbitrary, then one obtains as in section 3 for $t \ge 1$
$$
\align
\|w(\cdot,t)\|_{4,m} &\le Ct^{m-k} \epsilon(u_0) + \int_0^1 C (1+t-s)^{m-k} N_1(t-s)
\|w(\cdot,s)\|_{2,m} \, ds \\
& \quad + \int_1^t C(1+t-s)^{m-k} N_1(t-s) \|G(v)(\cdot,s)\|_{2,k} \, ds \\
&\le Ct^{m-k} \epsilon(u_0) + \int_1^t C (1+t-s)^{m-k} N_1(t-s) \|w\cdot,s)\|_{4,k/2}^2 \, ds \\
& \le Ct^{m-k} \epsilon(u_0) \, .
\endalign
$$
Finally for $t\ge 1$, 
$$
\align
\|w(\cdot,t)\|_{2,m} &\le Ct^{m-k} \epsilon(u_0) + \int_0^t C (1+t-s)^{m-k} N_0(t-s)
\|w(\cdot,s)\|_{4,k/2}^2 \, ds \\
&\le Ct^{m-k} \epsilon(u_0) + \int_0^t C (1+t-s)^{m-k} N_0(t-s) (1+s)^{-k} 
\epsilon^2(u_0) \, ds \\
& \le Ct^{m-k} \epsilon(u_0) \, .
\endalign
$$
The theorem is now completely proved. 
\enddemo

\medskip
Let us note in concluding that estimates  (B.7a,b) can easily be modified to 
$$
\|S(t)\phi\|_{\infty,k} \le C t^{-1/4}\|\phi\|_{2,k}, \quad
\|[S(t)\phi]_x\|_{\infty,k} \le C t^{-3/4}\|\phi\|_{2,k}
$$
for $0 < t \le 1$.  The last argument in the proof of Theorem 4.1 then implies that also
$$
\|u(\cdot,t) - \phi\|_{\infty,m} \le C t^{m-k} \epsilon(u_0) \, . 
$$

\bigskip
\head
Appendix A: Linearization of Scalar Viscous Conservation Laws
\endhead

In this section, properties of solutions of the equation
$$
u_t - u_{xx} + c u_x + du = 0 \quad (x \in \bold R, t > 0) \tag A.1
$$
are collected that are used in the main part of this paper.  
Here  $c$ and  $d$ are suitable coefficient functions which are bounded together with their first derivatives.  The results are mostly well known.

\medskip
Define the Banach spaces
$X = L^\infty(\bold R) \cap C^0(\bold R)$ and $X_0 = \{y \in X \, \big| \,
y \,\, \text{is uniformly continuous} \, \}$, both equipped with the
supremum norm $\| \cdot\|_\infty$. Define the operator $L$ acting in $X$
by $L\varphi = -\varphi_{xx} + c\varphi_x + d \varphi$ for $\varphi \in D(L) = \{ y \in X \,
\| \, y'' \in X \}$.  The restriction of $L$ to $X_0$ will also be denoted 
by $L$. Then $-L$ generates a $C_0$ semigroup $(S(t))_{t \ge 0}$ in 
$X_0$ which can be constructed as a perturbation of the heat semigroup.
In addition, the estimates hold 
$$
\align
\|S(t)\varphi\|_\infty &\le \|\varphi\|_\infty \qquad (0 \le t < \infty) 
\quad  \text{if} \quad d = 0 \tag A.2\\
\|\left[ S(t)\varphi\right]_x\|_\infty &\le \frac{C}{\sqrt{t}}\|\varphi\|_\infty
\qquad (0<t \le 1) \tag A.3
\endalign
$$
The first estimate is the maximum principle, the second follows from 
estimates for fundamental solutions in \cite{9}. The semigroup can be extended
to act on $X$, and the extension will also be denoted by $S(t)$.  
The extension still satisfies the estimates above. It is 
only strongly continuous for $t>0$, but this is irrelevant for the purposes of this paper.  

\medskip
Let us next examine the behavior of this semigroup in subspaces of $X_0$ that are defined by means of weight functions.  Let  $w: \bold R
\, \to \, [1,\infty)$ be smooth, with $w(\pm \infty) = \infty$, and with the 
first three derivatives of $x \, \mapsto \, \log \, w(x)$ bounded.
Set $X_w = \{y \in X \, \big | \, y \cdot w \in X \}$ with the norm
$\|y\|_{\infty,w} = \|w \cdot y \|_\infty$. There is a natural bijection
$R: X_w \, \to \, X, \, Ry = wy$, inducing the (formal) conjugate
$\tilde S(t) = RS(t)R^{-1}$ of $S(t)$ which again acts on $X$.  A straight forward
computation shows that the restriction of $\tilde S(t)$ to $X_0$ has an 
infinitesimal generator $-\tilde L$ of the same form as $L$.  
The estimates in \cite{9} apply to
this more general case and imply that $\tilde S(t)$ acts on 
$X_0$ as a $C_0$ - semigroup.  Thus $S(t)$ acts
on each $X_w$ with estimates $\|S(t) y \|_{\infty,w} \le C e^{Mt}\|y\|_{\infty,w}$ for 
$0 \le t < \infty$ and $\|\left[ S(t)y\right]_x \|_{\infty,w} \le \dfrac{C}
{\sqrt{t}} \|y\|_{\infty,w} $ for $0 < t \le 1$.  The constants $C$ and $M$ now 
depend also on $w$.  The same notation for $S(t)$ will be used, whether it acts 
on $X$ or on $X_w$.

\medskip
Now let $\rho, \, k > 0$ and consider specifically the spaces $X_\rho$ and $X_k$ 
that correspond to the weight functions $w(x) = \cosh(\rho x)$ and 
$w(x) = \left(1+x^2\right)^{k/2}$ (as in section 1). The work in \cite{4} and \cite{15}
implies that for sufficiently small $\rho$ there exists $\epsilon = 
\epsilon(\rho)>0$ and $C>0$ such that for all $\varphi \in X_\rho$ and all $t \ge 0$
$$
\|S(t)\varphi \|_{\infty,\rho} \le C e^{-\epsilon t}\|\varphi\|_{\infty,\rho} \,. \tag A.4
$$
Theorem 2.2 and (A.2) then imply that for all $\varphi \in X_k$ and all $0 < m < k$ 
$$
\|S(t)\varphi \|_{\infty,m} \le C_1(1+t)^{m-k}\|\varphi \|_{\infty,k} \,. \tag A.5
$$
Using (A.3), it also follows that for all $t>0$
$$
\|\left[S(t)\varphi \right]_x \|_{\infty,m} \le N(t) C_1(1+t)^{m-k}\|\varphi\|_{\infty,k}  \tag A.6
$$
with $N(t) = \max \{ 1, \, t^{-1/2} \}$. Finally, since $[S(t)\varphi]_x$ satisfies a parabolic equation of the same form as (A.1), the estimates hold for $0 < t \le 1$ 
$$
\|[S(t)\varphi]_x\|_{\infty,k} \le CN(t)\|[\varphi]_x\|_{\infty,k} \quad \text{and} \quad
\|[S(t)\varphi]_x\|_\infty \le C\|[\varphi]_x\|_\infty \, .
\tag A.7
$$

\bigskip

\head
Appendix B: Linearization of Generalized Korteweg-De Vries - Burgers Equations
\endhead

This appendix collects properties of solutions of the equation
$$
u_t + u_{xxx} - \alpha u_{xx} + c u_x + du = 0 \quad (x \in \bold R, t > 0) \tag B.1
$$
which are used in the main part of the paper.  Here  $\alpha > 0$, and $c$ and  $d$ are 
suitable smooth coefficient functions. Proofs will only be indicated.

\medskip
If $c = d = 0$, the spatial Fourier transform $\hat u(\cdot,t)$ of the solution $u$ of (B.1)
is given by
$$
\hat u(\xi,t) = e^{(i\xi^3- \alpha \xi^2)t} \hat u(\xi,0) \, \tag B.2
$$
This defines a contraction semigroup $S_0(t)_{t \ge 0}$ in $L^2$.  Moreover, for $\rho < 
\alpha/3$, $S_0$ maps $L^{2,\rho}$ into itself, since Fourier transforms of functions in 
$L^{2,\rho}$ have analytic extensions into the strip $\{z\, | \, |\Im(z)| < \rho \, \}$ 
that are square integrable on its boundary and since the multiplier in (B.2) is bounded by 
$C e^{-(\alpha - 3\rho)|\Re(\xi)|^2t}$ 
on any such strip.  Thus $S_0$, restricted to any such $L^{2,\rho}$, also generates a $C_0$
semigroup there.
In addition, since the multiplier  $\xi^2 t e^{(i\xi^3- \alpha\xi^2)t}$ is similarly
bounded, there are the estimates for all  $t$
$$
\|[S_0(t)\varphi]_{xx}\|_2 \le \frac{C}{t}\|\varphi\|_2 \quad \text{and} \quad
\|[S_0(t)\varphi]_{xx}\|_{2,\rho} \le \frac{C}{t}\|\varphi\|_{2,\rho} \tag B.3
$$
and from standard calculus estimates
$$
\|[S_0(t)\varphi]_{x}\|_2 \le \frac{C}{\sqrt{t}}\|\varphi\|_2
\quad \text{and} \quad
\|[S_0(t)\varphi]_{x}\|_{2,\rho} \le \frac{C}{\sqrt{t}}\|\varphi\|_{2,\rho}  \,.\tag B.4
$$
By Theorem 2.1, $S_0$ therefore maps each $L^{2,k}$ into itself and satisfies similar estimates for  $0 < k < \infty$.

\medskip
By a standard perturbation argument, one obtains existence and uniqueness of solutions of 
(B.1)for general coefficients $c,\, d$ that are bounded together with their first and second
derivatives.  The resulting semigroup $S(t)_{t \ge 0}$ maps $L^2, \,L^{2,\rho}\,  
(\rho < 1/3), \, L^{2,k} \, 
(0 < k < \infty)$ into itself, and the estimates hold
$$
\align
\|S(t)\varphi\|_{2,k} &\le C\|\varphi\|_{2,k} \tag B.5a \\
\|[S(t)\varphi]_{x}\|_{2,k} &\le \frac{C}{\sqrt{t}}\|\varphi\|_{2,k} \tag B.5b  \\
\|[S(t)\varphi]_{xx}\|_{2,k} &\le \frac{C}{t}\|\varphi\|_{2,k}  \, .\tag B.5c
\endalign 
$$

\medskip
Next note that for differentiable initial data $u(\cdot,0)$, the derivative $u_x$ 
also satisfies an equation of the form (B.1).  Thus there is also the estimate
$$
\|[S(t)\varphi]_{x}\|_2 \le C\| \varphi_x\|_2  \quad (0 < t \le 1) \tag B.6
$$
Finally, the estimates
$$
\align
\|S(t)\varphi\|_{4,k} \le Ct^{-1/8} \|\varphi\|_{2,k} \quad (0 < t \le 1) \quad \tag B.7a \\
\|[S(t)\varphi]_{x}\|_{4,k} \le Ct^{-5/8} \|\varphi\|_{2,k} \quad (0 < t \le 1) \quad \tag B.7b
\endalign
$$
for all $\varphi \in L^{2,k}$ follows from (B.5a,b,c) and the calculus inequality
$\|v\|^4_{4,k} \le C \|v\|^3_{2,k} \left( \|v\|_{2,k}+\|v_x\|_{2,k} \right) $.

\medskip
Let us next turn to estimates for $S(t)$ for large $t$, in the special case where $d = 0$.
Consider a general weight function $w > 0 $ with the properties
$$
|w'| \le \kappa w, \quad |w''| + |w'''| \le C, \quad \text{with} \quad \kappa < \alpha/3, \, 
C > 0 \, .
$$
A calculation shows that for $u(\cdot,t) = S(t)\varphi$ with $\varphi \in L^{2,\rho}, \, 
\rho < \alpha/3$, and $t>0$
$$
0 = \frac{d}{dt} \int w^2u + \int u_x^2 \left(3 |w'|^2 + 2\alpha w^2 \right) -
\int u^2 \left( \alpha (w^2)'' + (w^2)''' + (cw^2)' \right) \, . \tag B.8
$$
Consider first the case $w = 1$. Then this identity implies
$$
\frac{d}{dt} \int u^2 \le \int c' u^2 \, . \tag B.9
$$ 
Therefore if  $c' \le 0$, then 
$$
\|S(t)\varphi\|_2 \le \|\varphi\|_2  \tag B.10
$$
as was observed in \cite{14}. Next fix the assumptions
$$
c'(x)  < 0 \quad (x \in \bold R), \lim_{x \to -\infty} c(x) = \quad c_L  > 0 > 
c_R = \lim_{x \to +\infty} c(x) \, .  \tag B.11
$$
Under these assumptions, for all sufficiently small $\rho > 0$ there
are $C, \gamma > 0$ such that
$$
\|S(t)\varphi\|_{2,\rho} \le Ce^{-\gamma t}\|\varphi\|_{2,\rho} \, . \tag B.12
$$
Indeed, consider the weights $w(x) = \sqrt{\cosh(\rho(x - x_0))}$ where $x_0$ is such 
that $c(x_0) = 0$. From (B.8) one then obtains
$$
\frac{d}{dt} \int w^2 u^2 \le \int F_\rho w^2 u^2 \tag B.13
$$
where
$$
F_\rho(x) = \alpha \rho^2 + c'(x) + \left(\rho c(x) + \rho^3 \right) \tanh \rho (x-x_0) \, .
$$
Let us show that for all sufficiently small $\rho>0$, $F_\rho < -\gamma < 0$ on $\bold R$.
Set $2 A = \min (c_L, -c_R)$. If $M$ is sufficiently large (depending on $A$) and $\rho^2 < A$, 
then for $|x-x_0|>M$
$$
F_\rho(x) \le \alpha \rho^2 + (\rho^2 - A)\rho \tanh \rho M 
		\le - \gamma_0 \rho^2 
$$
for some $\gamma_0 > 0$.  If now $|x-x_0| \le M$, then
$$
F_\rho(x) \le \alpha \rho^2 + C \rho + \sup_{|x-x_0|\le M} c'(x) 
	     \le - \gamma_1
$$
provided $\rho$ is decreased further.  Then (B.13) implies
$\frac{d}{dt} \int w^2 u^2 \le -\gamma \int w^2 u^2 $, 
and (B.12) follows.  Theorem 2.1 and (B.5) now imply that for all $\varphi \in L^{2,k}$ and 
$0 < m < k, \, t > 0$
$$
\align
\|S(t)\varphi \|_{2,m} &\le C(1+t)^{m-k} \|\varphi\|_{2,k} \tag B.14a \\
\|[S(t)\varphi]_x \|_{2,m} &\le C (1+t)^{m-k} N_0(t) \|\varphi\|_{2,k} 
\tag B.14b
\endalign
$$
with  $N_0(t) = \max(1,t^{-1/2})$. Finally note that (B.7) and (B.13) imply
$$
\align
\|S(t)\phi\|_{4,m} &\le C t^{m-k}\|\phi\|_{2,k} \tag B.15a\\
\|[S(t)\phi]_x\|_{4,m} &\le C t^{m-k}\|\phi\|_{2,k} \tag B.15b
\endalign
$$
for $t \ge 1$, for all $0 < m < k$. 

\bigskip

\head
Appendix C: An Integral Estimate
\endhead

\proclaim{Lemma C.1}
Let $ 0 < \alpha < \beta $ with $ \beta >1$. Let 
$M:(0,\infty) \, \to \, \bold R$ be bounded on $[1, \infty)$ and 
integrable on $(0,1)$.
Then there exists a constant $C = C_{\alpha,\beta}$ such that for all
$t \ge 0$
$$
\int_0^t M(t-s) (1+t-s)^{-\alpha}(1+s)^{-\beta} \, ds \le Ct^{-\alpha} .
$$
\endproclaim

\demo{Proof}
Split the integral:
$$
\int_0^t M(t-s) (1+t-s)^{-\alpha}(1+s)^{-\beta} \, ds \le 
\int_0^t (1+t-s)^{-\alpha}(1+s)^{-\beta} \, ds +
\int_0^1 M(s) (1+s)^{-\alpha}(1+t-s)^{-\beta} \, ds .
$$
The second integral is clearly bounded by $C(1+t)^{-\beta}$.  
The first integral becomes
$$
\int_0^t (1+t-s)^{-\alpha}(1+s)^{-\beta} \, ds  
	 = \int_0^{t/2} (1+t-s)^{-\alpha}(1+s)^{-\beta} \, ds + 
	\int_{t/2}^t (1+t-s)^{-\alpha}(1+s)^{-\beta} \, ds  \, .
$$
If $\alpha, \,\beta > 1$, one estimates further
$$
...  \le (\beta - 1)^{-1} (1+t/2)^{-\alpha} + (\alpha - 1)(1+t/2)^{-\beta} 
\le C(1+t)^{-\alpha}\, . 
$$
If $\alpha < 1 < \beta$, the same estimate yields
$$
... \le (\beta - 1)^{-1} (1+t/2)^{-\alpha} + 
(1- \alpha)(1+t)^{1-\alpha}(1+t/2)^{-\beta} 
\le C(1+t)^{-\alpha} \, .
$$
Finally if $\alpha = 1 < \beta$, then the modified estimate holds
$$
...	 \le (\beta - 1)^{-1} (1+t/2)^{-1} + 
        \log(1+t/2)(1+t/2)^{-\beta}  \le C (1+t)^{-1}   \, .
$$
This proves the lemma. A closer look at the proof shows that the exponent $\alpha$ on the 
right hand side cannot be improved. 
\enddemo

\enddocument